\documentclass[12pt]{article}
\usepackage[margin=1in]{geometry}

\usepackage[T1]{fontenc}
\usepackage{amsthm, amsmath, amsfonts, amscd}
\usepackage[symbol]{footmisc}
\usepackage[shortlabels]{enumitem} 
\usepackage{witharrows}
\WithArrowsOptions{tikz={font=\small}, displaystyle}

\usepackage{authblk}
\makeatletter
\renewcommand{\maketitle}{
	\begin{center}
		{\Large\textsc\@title\par}
			\vskip 0.5em			
        {\normalsize\@author\par}
		{\small\itshape\@date\par}
		\vskip 0.5em		
	\end{center}
	\hrule	
	\vskip 2em
}
\makeatother

\usepackage{titlesec} 
\titleformat{\section}{\large\bfseries\scshape}{\thesection.}{0.5em}{}[{\vspace{0.5em}\titlerule \vspace{0.5em}}]
\titleformat*{\subsection}{\large\scshape}

\theoremstyle{plain}
\newtheorem{theorem}{Theorem}
\newtheorem{corollary}{Corollary}

\newtheorem{proposition}{Proposition}

\theoremstyle{definition}
\newtheorem{definition}{Definition}

\theoremstyle{remark}
\newtheorem{innercase}{Case}
\newenvironment{case}[1]  {\renewcommand\theinnercase{#1}\innercase}
  {\endinnercase}

\def\N{\mathbb N} 
\def\R{\mathbb R} 
\def\Z{\mathbb Z} 
\def\P{\mathbb P} 
\def\E{\mathbb E} 
\def\F{\mathbb F} 

\def\Had{P_{\mathrm{Had}}} 
\renewcommand{\epsilon}{\varepsilon}
\def\supp{\mathrm{supp}}

\def\poly{\mathrm{poly}}

\usepackage{hyperref}
\usepackage{cleveref}
\usepackage{url}

\begin{document}

\title{Integer points in dilates of polytopes}
\author[1]{Shubhangi Saraf\footnote[1]{Research partially supported by a Sloan research fellowship and an NSERC Discovery Grant. Email: shubhangi.saraf@utoronto.ca}} 
\author[2]{Narmada Varadarajan\footnote[2]{Email: narmada.varadarajan@mail.utoronto.ca}}
\affil[1]{Department of Computer Science, University of Toronto}
\affil[2]{Department of Mathematics, University of Toronto}
\date{}
\maketitle

\begin{abstract}

In this paper we study how the number of integer points in a polytope grows as we dilate the polytope. 
We prove new and essentially tight bounds on this quantity by specifically studying dilates of the Hadamard polytope. 

Our motivation for studying this quantity comes from the problem of understanding the maximal  number of monomials in a factor of a multivariate polynomial of $s$ monomials. 
A recent result by Bhargava, Saraf and Volkovich~\cite{BSV18} showed that if $f$ is an $n$-variate polynomial, where each variable has degree $d$, and $f$ has $s$ monomials,
then any factor of $f$ has at most $s^{O(d^2 \log n)}$ monomials. 
The key technical ingredient of their proof was to show that any polytope with $s$ vertices, where each vertex lies in $\{0,..,d\}^n$, can have at most $s^{O(d^2 \log n)}$ integer points.  
The precise dependence on $d$ of the number of integer points was left open. 
We show that this bound---particularly the dependence on $d$---is essentially tight by studying dilates of the Hadamard polytope and proving new lower bounds on the number of its integer points. 
\end{abstract}

\section{Introduction}

Integer polytopes play a starring role in several areas of mathematics, and have recently found very interesting applications in theoretical computer science.
In particular, the problem of counting integer points in polytopes lies at the intersection of combinatorics and convex geometry, with applications in algebraic geometry \cite{ag1, ag2}, number theory \cite{nt1, nt2}, optimization \cite{opt1,opt2}, cryptography \cite{cry1}, and algebraic complexity theory~\cite{BSV18}, just to name a few.
Perhaps the most celebrated result on this topic is the \emph{Ehrhart polynomial}: given an integer polytope $P \subseteq \R^n$, the number of integer points in its \emph{dilate}, $|d P \cap \Z^n|$, is a polynomial in $d$ \cite{ehrhart}. 
When $P$ is an $n$-dimensional polytope, the leading term of its Ehrhart polynomial is $\mathrm{vol}(P) d^n$, which is a good approximation to the number of integer points for large values of $d$.
However, the coefficients of the Ehrhart polynomial are notoriously difficult to compute (see, for instance, \cite{ehrhart2} for examples of polytopes with \textit{negative} Ehrhart coefficients), so we have little information about integer points when the dilation factor $d$ is much smaller than the dimension $n$.

\subsection*{Motivation from algebraic complexity theory}
Our interest in estimating $|dP \cap \Z^n|$ when $d$ is smaller than $n$ comes from the problem of \emph{sparse polynomial factorization} in complexity theory.
A multivariate polynomial $f \in \F[X_1, \dots, X_n]$ where each variable $X_i$ appears with degree $\leq d$ is \emph{sparse} if it only has $\poly(n)$ monomials.
Polynomial factorization is a fundamental question in computational algebra, and sparse polynomial factorization is a special case that has been studied for several decades now \cite{BSV18, sparse1, sparse2}.
The first randomized algorithm for sparse polynomial factorization exhibited a dependence on the sparsity of the factors of the polynomial \cite{sparse2}, which raised the natural question: do factors of sparse polynomials have to be sparse?

We know that this is not true in general: if $g$ is a factor of some polynomial $f$, there is no polynomial upper bound for the number of monomials of $g$ in terms of the number of monomials of $f$. 
In characteristic $0$, we have the counterexample $f = \prod_{i=1}^{n}(X_i^{d} - Y_i^{d})$ and $g = \prod_{i=1}^n (X_i^{d-1} + X_i^{d-2}Y_i + \dots + Y_i^{d-1})$.
Here, $f$ has $2^n$ monomials, but $g$ is a factor of $f$ with $d^n$ monomials. 
In characteristic $p$, this gap can be more stark. 
Consider $f = \sum_{i=1}^n X_i^p$ and $g = \left( \sum_{i=1}^n X_i \right)^{\lfloor p/2 \rfloor}$, so $f$ has $n$ monomials, but $g$ is a factor of $f$ with at least $n^{p/2}$ monomials.
In general, we do not know better than the trivial upper bound of $d^n$ for the number of monomials in a factor of a sparse polynomial. 
However, note that in the setting where we assume that the individual degree $d$ is an absolute constant and $n$ can be arbitrarily large, these two counterexamples no longer exhibit significant sparsity gaps.

The current best-known bound for sparse factorization in this setting is due to \cite{BSV18}, who show that if $f$ has $s$ monomials and $g$ is a factor of $f$, then $g$ has at most $s^{O(d^2 \log n)}$ monomials.
The key technique of their proof is to reduce the problem of counting monomials in $f$ and $g$ into one of counting integer points in their Newton polytopes. 
In this paper we further study the problem of counting integer points in polytopes, since a better upper bound there would directly imply a better bound for the factoring problem.
Unfortunately, our results indicate that there is no better upper bound for polytopes, so any attempt to improve the bound for sparsity will require a new approach.

\subsection*{Connections to studying dilates of polytopes}

The key technique to proving the upper bound of  $s^{O(d^2 \log n)}$ in \cite{BSV18} is to study integer points in polytopes. 
In particular if $f = g h$, then one can show that the Newton polytope of $f$ is the Minkowski sum of the Newton polytopes of $g$ and $h$. 
The high-level idea in the proof is to show the following: If $g$ has too many monomials, then the Newton polytope of $g$ must necessarily have many vertices. 
This can then be used to show that the Newton polytope of $f$ must also have many vertices, hence $f$ must have many monomials. 
Their main technical contribution is to show that any polytope with $s$ vertices, where each vertex lies in $[-d,d]^n$, can have at most $s^{O(d^2 \log n)}$ integer points.  
The precise dependence on $d$ of the number of integer points was left open.

The main goal of our paper is to understand if this bound for integer points in polytopes is tight.
Saptharishi observed that the $\log n$ in the exponent of the bound is necessary (see~\cite{BSV18} for a proof of this) by counting integer points in \emph{Hadamard polytopes}, a class of simplices with vertices in $\{-1,1\}^n$.
Specifically, he shows that the number of integer points in the Hadamard polytope is at least $n^{\Omega(\log n)}$. 
In this paper, we show that the dependence of $d^2$ in the exponent is also essentially tight by studying dilates of the Hadamard polytope. 

An astute reader may recognize the Hadamard polytope as the extremal solution to Hadamard's maximal determinant problem: the parallelepiped spanned by its vertices maximizes volume among all parallelepipeds with vertices in $\{-1,1\}^n$.
We expect that the Hadamard polytope should also maximize the number of integer points---a discrete analog of volume---among all simplices with vertices in $\{-1,1\}^n$.
This is the heuristic that motivates us to study the Hadamard polytope as the ``worst-case simplex'' for our problem.

In the case when $s = \poly (n)$ (as with factors of sparse polynomials), we can reduce our problem about counting integer points in polytopes in $[-d,d]^n$ to counting integer points in dilates of polytopes as follows.
First, by Carath\'eodory's theorem, each vertex of $P$ is a convex combination of at most $n+1$ points of $\{-d,d\}^n$.
Replacing the vertices of $P$ with these points, we obtain a polytope $P_d$ containing $P$ that still has $\poly(n)$ vertices, but they now lie in $\{-d,d\}^n$.
So, $P_0 = \frac{1}{d}P_d$ is an integer polytope with vertices in $\{-1,1\}^n$ and $|dP_0 \cap \Z^n| \geq |P \cap \Z^n|$.

Our contribution in this paper is to provide lower bounds for the number of integer points in the Hadamard polytope.
Our first main result is an elegant algebraic characterization of the integer points in the Hadamard polytope by constructing an explicit bijection with affine subspaces of $\F_2^{\log n}$.
This underpins part of our second main result: a lower bound for the number of integer points in small dilates.
To lower bound the points in larger dilates, we use the probabilistic method.
These bounds show that both the dependences on $d^2$ and on $\log n$ in the bound $|d P \cap \Z^n| \leq |V(P)|^{O(d^2 \log n)}$ are necessary.

\section{Main results}

Let $n = 2^m$ for some $m \in \N$.
The $n \times n$ \emph{Hadamard matrix} $H$ has its rows and columns indexed by elements of the vector space $\F_2^m$:
\begin{align*}
    H(a,b) &= (-1)^{\langle a,b \rangle},
\end{align*}
where $\langle a,b \rangle$ is the standard dot product on $\F_2^m$.

\begin{definition}
    The \emph{Hadamard polytope} $\Had \subset \R^n$ is the convex hull of the column vectors of $H$.
\end{definition}

The columns of $H$ are linearly independent so they form the vertices of $\Had$; the \emph{Hadamard polytope} is an $(n-1)$-dimensional simplex with vertices in $\{-1, 1\}^n$.
The lower bound $|\Had \cap \Z^n| \geq n^{\Omega(\log n)}$ was first observed by Ramprasad Saptharishi and appears in the paper \cite{BSV18}.
Our contribution is a matching upper bound on $|\Had \cap \Z^n|$ through a complete characterization of all integer points in the Hadamard polytope using the vector space structure of $\F_2^m$.
This characterization will also be useful when we study integer points in the dilates of $\Had$.

First, some notation: let $\{h_a : a \in \F_2^m\}$ be the set of column vectors of the Hadamard matrix $H$, i.e.\ the vertices of the Hadamard polytope $\Had$.
We can also identify the coordinates of $\R^n$ with vectors in $\F_2^m$, with the convention that the first coordinate in $\R^n$ corresponds to the zero vector in $\F_2^m$.
For every $a \in \F_2^m$ and $v \in \R^n$, let $v(a)$ denote the $a$th coordinate of $v$ in the standard basis, and 
\begin{align*}
   \supp(v) = \{a \in \F_2^m: v(a) \neq 0\} 
\end{align*}
We also know that every point $v \in \Had$ can be expressed as a convex combination $v = \sum_{a \in \F_2^m}t_a h_a$.
Since the vectors $\{h_a : a \in \F_2^m\}$ are linearly independent, this expression is unique, and we can define
\begin{align*}
    T(v) &= \{a \in \F_2^m: t_a \neq 0\}.
\end{align*}
This can also be thought of as the support of $v$ in the basis $\{h_a: a \in \F_2^m\}$.

\subsection*{Statements of main results}
\begin{proposition}\label{prop1}
    For every $v \in \Had \cap \Z^n$, $\supp(v)$ is a subspace of $\F_2^m$.
    Further, $T(v) = \supp(v)^{\perp} + b$ for some $b \in \F_2^m$, and 
    \begin{align*}
        v &= \frac{1}{|T(v)|}\sum_{c \in T(v)}h_c
    \end{align*}
    is actually a uniform linear combination of vertices.
\end{proposition}
\begin{corollary}\label{cor1}
    The number of integer points in $\Had$ is $|\Had \cap \Z^n| = n^{\Theta(\log n)}$.
\end{corollary}
\begin{theorem}\label{thm1}
    For any $d < n \log n$ and fixed $0 < \epsilon < 1/2$, the number of integer points in the dilate $d\Had$ is lower bounded by
    \begin{align*}
        |d\Had \cap \Z^n| 
        \geq \begin{cases}
        n^{\Omega(d \log n)}, & d \leq \frac{\log n}{4};\\
        n^{\Omega(\log^2n)}, & \frac{\log n}{4} \leq d \leq \log^{1.5}n\\
        n^{\epsilon d^2/2\log n}, & \log^{1.5}n \leq d \leq (n\log n)^{1/2 - \epsilon};\\
        \left( d^{2\epsilon}/4 \right)^n, & (n\log n)^{1/2 + \epsilon} \leq d < n \log n.
        \end{cases}     
    \end{align*}
\end{theorem}

\subsection*{Proofs of main results}
\begin{proof}[Proof of \Cref{prop1}]
    First, we need to show that $\supp(v)$ is a subspace of $\F_2^m$.
    We will actually show something stronger: if $v(a) \neq 0$ and $v(b) \neq 0$, then $v(a+b) = v(a) v(b)$.
    Note that since $v$ is an integer point in $\Had \subseteq [-1,1]^n$, its coordinates actually lie in $\{-1,0,1\}^n$.
    So, the three cases we need to consider for our proof are: first, $v(a) = v(b) = 1$; second, $v(a) = v(b) = -1$; and third, $v(a) = 1$, $v(b) = -1$.
    We will write out the details for the first case when $v(a) = v(b)=1$, and the other cases will follow a similar computation.
    
    Since $v$ is a convex combination of the vertices of $\Had$, we can write it as $v = \sum_{c \in \F_2^m}t_ch_c$, for some nonnegative real numbers $(t_c)_{c \in \F_2^m}$ that satisfy $\sum_{c \in \F_2^m}t_c = 1$.
    From the definition of the Hadamard matrix $H$, we get the following identities for the coordinates $a,b \in \F_2^m$,
    \begin{align*}
        v(a) &= \sum_{c \in \F_2^m}t_ch_c(a) = \sum_{c \in \F_2^m}(-1)^{\langle c,a\rangle}t_c = \sum_{c \perp a}t_c - \sum_{c \not\perp a}t_c,\\
        v(b) &= \sum_{c \in \F_2^m}t_ch_c(b) = \sum_{c \in \F_2^m}(-1)^{\langle c,b\rangle}t_c =  \sum_{c \perp b}t_c - \sum_{c \not\perp b}t_c.\\
    \end{align*}
    
    Define
    \begin{align*}
        s_{ab}= \sum_{c \perp a, c \perp b} t_c;\quad 
        s_{a}= \sum_{c \perp a, c \not\perp b} t_c;\quad 
        s_{b}= \sum_{c \not\perp a, c \perp b} t_c;\quad 
        s_{0}= \sum_{c \not\perp a, c \not\perp b} t_c.
    \end{align*}

    Along with the identity $\sum_{c \in \F_2^m} t_c=1$, this gives us the system of linear equations 
    \begin{align*}
     v(a)=s_{ab}+s_a-s_b-s_0 &=1\\
    v(b)=  s_{ab} -s_a + s_b - s_0 &=1\\
      s_{ab} + s_a + s_b + s_0 &=1\\
    \end{align*}
    The general solution to this system has the form
    \begin{align*}
        s_{ab} - s_0 &=1\\
        s_a + s_b &= 0\\
        s_b + s_0 &= 0.
    \end{align*}
    Since each of $s_{ab}, s_a, s_b,$ and $s_0$ must be nonnegative, this yields the unique solution $s_{ab}=1, s_a=s_b=s_0=0$.
    So,
    \begin{align*}
        v(a+b) &= \sum_{c \in \F_2^m}(-1)^{\langle a+b,c \rangle}t_c\\
        &= \sum_{c \perp a, c \perp b} t_c - \sum_{c \perp a, c \not\perp b} t_c - \sum_{c \not\perp a, c \perp b} t_c + \sum_{c \not\perp a, c \not\perp b} t_c\\
        &=s_{ab} - s_a - s_b + s_0 \\
        &= 1
        = v(a)v(b).
    \end{align*}
    The argument for the other cases follows similarly: we set up a system of equations in the variables $s_{ab}, s_b, s_a, s_0$ that has a unique solution subject to nonnegativity.
    Computing the value of $v(a+b)$ for this unique solution shows that $v(a+b) = v(a)v(b)$.
    This shows not only that $\supp(v)$ is a subspace of $\F_2^m$, but also that the entries of $v$ are determined by their values on a basis of $\supp(v)$.

    Now, we need some linear algebra to complete the proof.
    Let $b_1, \dots, b_k \in \F_2^m$ be a basis for $\supp(v)$, and choose $\epsilon_i \in \{0,1\}$ so that $v(b_i) = (-1)^{\epsilon_i}$ for each $i=1, \dots, k$.
    Consider the linear map $\F_2^m \to \supp(v)$ that sends $c \in \F_2^m$ to $\langle b_1,c \rangle b_1 + \dots + \langle b_k,c\rangle b_k \in \supp(v)$.
    The kernel of this map is exactly $\supp(v)^{\perp}$, so its range must be all of $\supp(v)$.
    In particular, there is some $b \in \F_2^m$ whose image under this map is $\epsilon_1 b_1 + \dots + \epsilon_k b_k$.
    So, $v(b_i) = (-1)^{\langle b_i,b \rangle}$ for each $i= 1,\dots,k$.
    For any $a \in \supp(v)$, we know that $a = \sum_{i \in I}b_i$ for some $I \subseteq [k]$, so
    \begin{align*}
        v(a) = \prod_{i \in I}v(b_i) = (-1)^{\sum_{i \in I}\langle b_i,b \rangle} = (-1)^{\langle a,b \rangle}.
    \end{align*}
    We claim that $T(v) = \supp(v)^{\perp} + b$ for this vector $b$.

    Define the vector $v' \in \R^n$ by
    \begin{align*}
        v' &= \frac{1}{\left|\supp(v)^{\perp}\right|}\sum_{c \in \supp(v)^{\perp} + b} h_c.
    \end{align*}
    We will show that $v' = v$ by showing that $v'(a) = v(a)$ for every coordinate $a \in \F_2^m$.
    
    First, for each $a \in \supp(v)$ and $c \in \supp(v)^{\perp} + b$, we have that $\langle a,c \rangle = \langle a,b \rangle$.
    So,
    \begin{align*}
        v'(a) = \frac{1}{|\supp(v)^{\perp}|}\sum_{c \in \supp(v)^{\perp} + b} h_c(a) = \frac{1}{|\supp(v)^{\perp}|}\sum_{c \in \supp(v)^{\perp} + b} (-1)^{\langle a,c \rangle} =(-1)^{\langle a,b \rangle} = v(a).
    \end{align*}
    Next, for each $a \notin \supp(v)$, the $\F_2$-linear map from $\supp(v)^{\perp} \to \{0,1\}$ that sends $c \in \supp(v)^{\perp}$ to $\langle a,c \rangle$ is not identically zero. 
    So, the size of the preimage of $0$ under this map is equal to the size of the preimage of $1$.
    This property still holds when we consider the translate $\supp(v)^{\perp} + b$: the dot product $\langle a,c\rangle$ is $0$ for exactly half the vectors $c \in \supp(v)^{\perp} + b$, and is $1$ for the other half.
    Thus, the average of $h_c(a)=(-1)^{\langle a,c \rangle}$ over all $c \in \supp(v)^{\perp} + b$ is zero, so
    \begin{align*}
        v'(a) = \frac{1}{|\supp(v)^{\perp}|} \sum_{c \in \supp(v)^{\perp}+b}(-1)^{\langle a,c \rangle} = 0 = v(a),
    \end{align*}
    since $a \notin \supp(v)$.
    This shows that $v$ does indeed have the form we claimed, with $T(v) = \supp(v)^{\perp} + b$ and
    \begin{align*}
        v = \frac{1}{|T(v)|}\sum_{c \in T(v)}h_c
    \end{align*}
        \end{proof}

\begin{proof}[Proof of \Cref{cor1}]
    Note again that the lower bound $|\Had \cap \Z^n| = n^{\Omega(\log n)}$ was already known by Saptharishi's contribution to \cite{BSV18}; we will provide a matching upper bound of $n^{O(\log n)}$.
   By \Cref{prop1}, integer points of $\Had \cap \Z^n$ are in one-to-one correspondence with affine subspaces of $\F_2^m$.
   Any $k$-dimensional linear subspace $V \subseteq \F_2^m$ has exactly $2^{m-k}$ corresponding affine subspaces, and there are exactly $\genfrac{[}{]}{0pt}{}{m}{k}_2$ (the Gaussian binomial coefficient) linear subspaces of dimension $k$.
   So, using the upper bound $\genfrac{[}{]}{0pt}{}{m}{k}_2 \leq 2^{4m^2}$ (see, for example, \cite{bin} for a proof of this bound),
   \begin{DispWithArrows*}
       |\Had \cap \Z^n| = \sum_{k=0}^m 2^{m-k}\genfrac{[}{]}{0pt}{}{m}{k}_2 \leq (m+1)2^{m+4m^2} = n^{O(\log n)},
   \end{DispWithArrows*}
   which completes the proof.
\end{proof}

Now we are ready to prove our main theorem: lower bounding the growth of integer points in the dilates $d \Had$.

\begin{proof}[Proof of \Cref{thm1}]
\begin{case}{1}
    $d \leq (\log n)/4$.
\end{case}
Let us introduce some more notation for our counting argument.
We will show that there are $n^{\Omega(d \log n)}$ distinct sets of $d$ integer points in $\Had$ whose sums are all distinct in $d\Had$.
For $v \in \Had \cap \Z^n$, define the \emph{exponent} of $v$ as $\exp(v) = \dim(\supp(v)^{\perp})$, so $|T(v)| = 2^{\exp(v)}$, and
\begin{align*}
    v &= \frac{1}{2^{\exp(v)}}\sum_{a \in T(v)}h_a.
\end{align*}
Let $t_a(v)$ be $2^{-\exp(v)}$ if $a \in T(v)$, and $0$ otherwise; that is, $t_a(v)$ is the coefficient of the vertex $h_a$ in the convex combination for $v$. 
Look at all sets $\{v_1, \dots, v_d\} \subseteq \Had \cap \Z^n$ such that
\begin{enumerate}[(i)]
    \item $T(v_i)$ is a linear subspace, i.e.\ $T(v_i) = \supp(v_i)^{\perp}$;
    \item $\frac{\log n}{2} - d < \exp(v_i) < \frac{\log n}{2} + d$;
    \item $\exp(v_i) \neq \exp(v_j)$ if $i \neq j$.
\end{enumerate}
Since $d \leq \frac{\log n}{4}$, for $n$ large enough, the number of subspaces of $\F_2^{\log n}$ with dimension $k$ for some $k \in \left(\frac{\log n}{2}-d, \frac{\log n}{2} + d\right)$ is $\genfrac{[}{]}{0pt}{}{\log n}{k}_2 \geq 2^{k(\log n - k)}= n^{\Omega(\log n)}$. 
So, the number of choices for the set of subspaces $\{T(v_1), \dots, T(v_d)\}$ with distinct exponents in $\left(\frac{\log n}{2}-d, \frac{\log n}{2} + d\right)$ is $n^{\Omega(d \log n)}$. 
We will now show that each set $\{v_1, \dots, v_d\}$ satisfying these properties produces a distinct point in $d\Had \cap \Z^n$, proving the lower bound $|d\Had \cap \Z^n| \geq n^{\Omega(d \log n)}.$

Now, suppose $\{v_1, \dots, v_d\}$ and $\{u_1, \dots, u_d\}$ are two distinct sets satisfying the above conditions.
And suppose for contradiction that $\sum_i v_i = \sum_i u_i$.
Expanding this in terms of the vertices of $\Had$, we get
\begin{align*}
    \sum_{i =1}^d\sum_{a \in \F_2^m}t_a(v_i)h_a &= \sum_{i =1}^d\sum_{a \in \F_2^m}t_a(u_i)h_a
\end{align*}
Since the vertices of $\Had$ are linearly independent, we must have the following stronger identity for all $a \in \F_2^m$:
\begin{align*}
    \sum_{i=1}^dt_a(v_i) = \sum_{i=1}^dt_a(u_i).
\end{align*}
For each fixed $a$, since $t_a(v_i)$ is either zero or $2^{-\exp(v_i)}$ and the exponents are all distinct, we can think of each sum as the unique binary representation of some dyadic rational number.

When $a =0$, we get $t_0(v_i) = 2^{-\exp(v_i)}$ and $t_0(u_i) = 2^{-\exp(u_i)}$ since $T(v_i)$ and $T(u_i)$ are linear subspaces containing $0$, for every $i=1,\dots, d$. 
This gives us the identity
\begin{align*}
    \sum_{i=1}^d 2^{-\exp(v_i)} = \sum_{i=1}^{d}2^{-\exp(u_i)}.
\end{align*}
So, the sets $\{ \exp(v_i): i \leq d \}$ and $\{ \exp(u_i): i \leq d \}$ must be equal since the above sums represent the same binary number.
This means that the only way for the sums $\sum_{i=1}^d v_i$ and $\sum_{i=1}^d u_i$ to be equal at the $0$th coordinate is if their exponents are all the same.
Now we will find another coordinate where they differ.

Assume without loss of generality that $\exp(v_i) = \exp(u_i)$ for all $i=1,\dots,d$, and $\exp(v_1) < \dots < \exp(v_d)$.
There must be some $k$ such that $T(v_k) \neq T(u_k)$.
Choose $a \in T(v_k) \setminus T(u_k)$ (which is possible since they are different subspaces of the same dimension), and now look at the identity
\begin{align*}
    \sum_{i=1}^{d}t_a(v_i) = \sum_{i=1}^{d}t_a(u_i).
\end{align*}
Again, each sum is the unique binary representation of some dyadic rational number, but $t_a(v_k) \neq 0$ while $t_a(u_k) = 0$, so they cannot represent the same number, which is the desired contradiction.
    \begin{case}{2}
      $(\log n)/4 \leq d \leq \log^{1.5}n$ .
    \end{case}
    Our lower bound from case 1 is monotone in $d$ and hence we can just reuse it for $d \geq \frac{\log n}{4}$,
    \begin{align*}
        |d\Had \cap \Z^n| \geq  \left| \frac{\log n}{4}\Had\cap \Z^n \right| \geq n^{\Omega(\log^2 n)}.
    \end{align*}
\begin{case}{3}
    $d \geq \log^{1.5}n$. 
\end{case}
For this final case we will use the probabilistic method.
    Take the projection $\R^n \to \R^{n-1}$ that deletes the first coordinate of every vector, and let $P_0$ be the image of $\Had$ under this map.
    Recall that $\Had  \subseteq \R^n$ is an $(n-1)$-dimensional simplex whose vertices all have first coordinate equal to $1$, since the first coordinate corresponds to the zero vector in $\F_2^m$. 
    So, this projection is a bijection on $\Had$.
    In particular, integer points of $dP_0$ are in one-to-one correspondence with integer points of $d\Had$, so $|d\Had \cap \Z^n| = |dP_0 \cap \Z^{n-1}|$.
    
    The new polytope $P_0$ is now a full-dimensional simplex in $[-1,1]^{n-1}$. 
    Since the vertices of the original Hadamard polytope $\Had \subseteq \R^n$ were orthogonal, for each vertex $v \in V(P_0) \subseteq \R^{n-1}$, the face of $P_0$ that does not contain $v$ is defined by the hyperplane $\{ x \in \R^{n-1} : \langle x,v \rangle = -1\}$.
    In particular, we can write the simplex $P_0$ as an intersection of half-spaces:
    \begin{align*}
        P_0 &= \bigcap_{v \in V(P_0)}\{x \in \R^{n-1}: \langle x, v \rangle \geq -1\},
    \end{align*}
    and
     \begin{align*}
        dP_0 &= \bigcap_{v \in V(P_0)}\{x \in \R^{n-1}: \langle x, v \rangle \geq -d\}.
    \end{align*}
    Fix some set of $c$ coordinates $S \subset [n]$ and some $D \in [d]$ (with $c$ and $D$ to be determined later) and look at the discrete hypercube $H_S[D]$ supported on those coordinates,
    \begin{align*}
        H_S[D] = \left([-D,D]\setminus \{0\}\right)^{S} \times \{0\}^{[n] \setminus S}.
    \end{align*}
    We will use the probabilistic method to show that $P_0$ contains at least $1/2$ the points of $H_S[D]$ for every choice of $S$.
    Note that the hypercubes $H_S[D]$ are disjoint by construction, so 
    \begin{align*}
      |dP_0 \cap \Z^n| \geq \sum_{S \subseteq [n], |S| =c}\left|dP_0 \cap H_S[D]\right|.  
    \end{align*}
  Fix a vertex $v \in V(P_0)$, and a hypercube $H_S[D]$.
  Let $X$ be a uniform random variable taking values in $H_S[D]$.
  Since $v$ is a $\pm 1$-valued vector, $\E|\langle X,v \rangle| = 0$. 
  Hoeffding's inequality implies that
    \begin{align*}
        \P(|\langle X,v\rangle| > d) &\leq 2e^{-d^2/2cD}.
    \end{align*}
    By the union bound over all vertices $v \in V(P_0)$,
    \begin{align*}
        \P(X \notin dP_0) \leq  \sum_{v \in V(P_0)}\P (|\langle X,v\rangle|>d) \leq 2ne^{-d^2/2cD}.
    \end{align*}
  Assuming that $2cD \leq \frac{d^2}{2\log n}$, the above probability is at most $1/2$, so at least $1/2$ the points of $H_S[D]$ must be contained in $dP_0$. 
 Now, summing over the disjoint sets $H_S[D]$ for all $\binom{n}{c}$ choices of $S$,
    \begin{align*}
    \left|d\Had \cap \Z^n\right| = \left|dP_0 \cap \Z^{n-1}\right| \geq
        \sum_{S \subseteq [n], |S|=c}|H_S[D] \cap dP_0| \geq \frac{1}{2}\binom{n}{c}(2D)^c.
    \end{align*}
 We divide our analysis into two subcases.

\begin{case}{3a}
$d \leq \left( n \log n \right)^{\frac{1}{2} - \epsilon}$, which implies that $\frac{d^2}{\log n} \leq n^{1-2\epsilon}$
\end{case}
In this case, set $D=1$ and $c = \frac{d^2}{4 \log n}$, so 
\begin{DispWithArrows*}[wrap-lines]
    |d\Had \cap \Z^n| &\geq \frac{1}{2}\binom{n}{c}(2D)^c \Arrow{Approximating the binomial coefficient}\\
    &\geq \left( \frac{n}{c} \right)^c \Arrow{Since $c = \frac{d^2}{4\log n} \leq n^{1-2\epsilon}$}\\
    &\geq \left( \frac{n}{n^{1-2\epsilon}} \right)^{d^2/4\log n}\\
    &\geq n^{ \epsilon d^2/2\log n }
\end{DispWithArrows*}

\begin{case}{3b} 
$d \geq \left( n \log n \right)^{\frac{1}{2} + \epsilon}$.
\end{case}
For our probabilistic argument to work, we need $2cD \leq \frac{d^2}{2\log n}$, which rearranges to $D \leq \frac{d^2}{4c \log n}$.
Set $D = d^{2 \epsilon}/4$ and $c = n$.
Now,
\begin{DispWithArrows*}[wrap-lines]
    \frac{d^2}{4c \log n} &= \frac{d^2}{4n \log n}\Arrow{Subsituting $n \log n \leq d^{2/(1+2 \epsilon)}$}\\
    &\geq \frac{d^2}{4d^{2/(1+2\epsilon)}}\\
    &= \frac{d^{4\epsilon/(1+2\epsilon)}}{4}\Arrow{Since $\epsilon \leq 1/2$}\\
    &\geq \frac{d^{2\epsilon}}{4} = D
\end{DispWithArrows*}
as desired.
Now, 
\begin{align*}
    |d\Had \cap \Z^n| &\geq \frac{1}{2}\binom{n}{c}(2D)^c \geq \left( d^{2\epsilon}/4 \right)^n.
\end{align*}

\end{proof} 

\bibliographystyle{alphaurl}
\bibliography{bib}

\newcommand{\etalchar}[1]{$^{#1}$}
\begin{thebibliography}{HHTY19}

\bibitem[AH10]{opt1}
Iskander Aliev and Martin Henk.
\newblock Feasibility of integer knapsacks.
\newblock {\em SIAM Journal on Optimization}, 20(6):2978--2993, 2010.
\newblock \href {https://doi.org/10.1137/090778043} {\path{doi:10.1137/090778043}}.

\bibitem[AJSe24]{opt2}
Divesh Aggarwal, Antoine Joux, Miklos Santha, and Karol~W\k egrzycki.
\newblock Polynomial time algorithms for integer programming and unbounded subset sum in the total regime, 2024.
\newblock Submitted July 7, 2024; revised July 11, 2024.
\newblock \href {https://arxiv.org/abs/2407.05435} {\path{arXiv:2407.05435}}, \href {https://doi.org/10.48550/arXiv.2407.05435} {\path{doi:10.48550/arXiv.2407.05435}}.

\bibitem[BR07]{nt1}
Matthias Beck and Sinai Robins.
\newblock {\em Computing the Continuous Discretely}.
\newblock Undergraduate Texts in Mathematics. Springer, New York, NY, 1 edition, 2007.
\newblock \href {https://doi.org/10.1007/978-0-387-46112-0} {\path{doi:10.1007/978-0-387-46112-0}}.

\bibitem[BSV20]{BSV18}
Vishwas Bhargava, Shubhangi Saraf, and Ilya Volkovich.
\newblock Deterministic factorization of sparse polynomials with bounded individual degree.
\newblock {\em J. ACM}, 67(2), May 2020.
\newblock \href {https://doi.org/10.1145/3365667} {\path{doi:10.1145/3365667}}.

\bibitem[BV25]{sparse1}
P.~Bisht and I.~Volkovich.
\newblock On solving sparse polynomial factorization related problems.
\newblock {\em Computational Complexity}, 34(7), 2025.
\newblock \href {https://doi.org/10.1007/s00037-025-00268-5} {\path{doi:10.1007/s00037-025-00268-5}}.

\bibitem[CLS12]{nt2}
Sheng Chen, Nan Li, and Steven~V. Sam.
\newblock Generalized ehrhart polynomials.
\newblock {\em Transactions of the American Mathematical Society}, 364(1):551--569, 2012.
\newblock \href {https://doi.org/10.1090/S0002-9947-2011-05494-2} {\path{doi:10.1090/S0002-9947-2011-05494-2}}.

\bibitem[Cox96]{ag1}
David~A. Cox.
\newblock Recent developments in toric geometry, 1996.
\newblock URL: \url{https://arxiv.org/abs/alg-geom/9606016}, \href {https://arxiv.org/abs/alg-geom/9606016} {\path{arXiv:alg-geom/9606016}}.

\bibitem[CPS24]{ag2}
Luis Crespo, \'Alvaro Pelayo, and Francisco Santos.
\newblock Ewald's conjecture and integer points in algebraic and symplectic toric geometry, 2024.
\newblock URL: \url{https://arxiv.org/abs/2310.10366}, \href {https://arxiv.org/abs/2310.10366} {\path{arXiv:2310.10366}}.

\bibitem[Ehr62]{ehrhart}
Eugène Ehrhart.
\newblock Sur les polyèdres rationnels homothétiques à n dimensions.
\newblock {\em Comptes rendus de l'Académie des Sciences}, 254:616--618, 1962.

\bibitem[FLC{\etalchar{+}}25]{cry1}
Yansong Feng, Hengyi Luo, Qiyuan Chen, Abderrahmane Nitaj, and Yanbin Pan.
\newblock Computing asymptotic bounds for small roots in coppersmith's method via sumset theory.
\newblock In Yael Tauman~Kalai and Seny~F. Kamara, editors, {\em Advances in Cryptology -- CRYPTO 2025}, pages 3--32, Cham, 2025. Springer Nature Switzerland.

\bibitem[HHTY19]{ehrhart2}
Takayuki Hibi, Akihiro Higashitani, Akiyoshi Tsuchiya, and Koutarou Yoshida.
\newblock Ehrhart polynomials with negative coefficients.
\newblock {\em Graphs and Combinatorics}, 35(1):363--371, 2019.

\bibitem[NP95]{bin}
Peter~M. Neumann and Cheryl~E. Praeger.
\newblock Cyclic matrices over finite fields.
\newblock {\em Journal of the London Mathematical Society, Series 2}, 52(2):263--284, 1995.
\newblock \href {https://doi.org/10.1112/jlms/52.2.263} {\path{doi:10.1112/jlms/52.2.263}}.

\bibitem[vK85]{sparse2}
Joachim {von zur Gathen} and Erich Kaltofen.
\newblock Factoring sparse multivariate polynomials.
\newblock {\em Journal of Computer and System Sciences}, 31(2):265--287, 1985.
\newblock URL: \url{https://www.sciencedirect.com/science/article/pii/0022000085900443}, \href {https://doi.org/10.1016/0022-0000(85)90044-3} {\path{doi:10.1016/0022-0000(85)90044-3}}.

\end{thebibliography}
\end{document}